\let\myfrac=\frac
\input eplain
\let\frac=\myfrac
\input amstex
\input epsf




\loadeufm \loadmsam \loadmsbm
\message{symbol names}\UseAMSsymbols\message{,}

\font\myfontdefault=cmr10

\font\mytdmchapfont=cmb10 at 14pt
\font\mytdmheadfont=cmb10 at 10pt
\font\mytdmsubheadfont=cmr10

\magnification 1200
\newif\ifinappendices
\newif\ifundefinedreferences
\newif\ifchangedreferences
\newif\ifloadreferences
\newif\ifmakebiblio
\newif\ifmaketdm

\undefinedreferencesfalse
\changedreferencesfalse


\loadreferencestrue
\makebibliofalse
\maketdmfalse

\def\headpenalty{-400}     
\def\proclaimpenalty{-200} 

%
%

\def\alphanum#1{\ifcase #1 _\or A\or B\or C\or D\or E\or F\or G\or H\or I\or J\or K\or L\or M\or N\or O\or P\or Q\or R\or S\or T\or U\or V\or W\or X\or Y\or Z\fi}
\def\gobbleeight#1#2#3#4#5#6#7#8{}

\newwrite\references
\newwrite\tdm
\newwrite\biblio

\newcount\chapno
\newcount\headno
\newcount\subheadno
\newcount\procno
\newcount\figno
\newcount\citationno

\def\setcatcodes{%
\catcode`\!=0 \catcode`\\=11}%

\ifloadreferences
    {\catcode`\@=11 \catcode`\_=11%
    \global\def\_@citation@BallSchGrom{1}
\global\def\_@citation@Brezis{2}
\global\def\_@citation@CaffNirSpruckV{3}
\global\def\_@citation@Gauss{4}
\global\def\_@citation@Gromov{5}
\global\def\_@citation@GuanSpruck{6}
\global\def\_@citation@ImaTani{7}
\global\def\_@citation@LabourieB{8}
\global\def\_@citation@Labourie{9}
\global\def\_@citation@LehtoVirtanen{10}
\global\def\_@citation@Martinez{11}
\global\def\_@citation@Pogorelov{12}
\global\def\_@citation@RosSpruck{13}
\global\def\_@citation@Sasaki{14}
\global\def\_@citation@SmiPKS{15}
\global\def\_@citation@SmiFCS{16}
\global\def\_@citation@SmiPPH{17}
\global\def\_@citation@TrudWang{18}
\global\def\_@proc@ThmCompactness{1.3}
\global\def\_@proc@ThmExistence{1.4}
\global\def\_@proc@ThmContinuousDependance{1.5}
\global\def\_@proc@LemmaQuasiConformal{2.1}
\global\def\_@proc@DefnCheegerGromov{3.2}
\global\def\_@proc@ThmLabCompactness{3.3}
\global\def\_@proc@LemmaSchwartz{3.4}
\global\def\_@proc@LemmaConvergenceMappingsConverge{3.5}
\global\def\_@proc@PropConformalRepresentation{4.1}
\global\def\_@proc@PropQuasiConfCompactness{4.2}
\global\def\_@proc@PropInversesConverge{4.3}
\global\def\_@proc@PropAvoidsLimitPoint{5.1}
\global\def\_@proc@PropStaysInCompactSet{5.2}
\global\def\_@proc@PropStrongCompactness{5.3}
    }%
\else
    \openout\references=references.tex
\fi

\newcount\newchapflag 
\newcount\showpagenumflag 

\global\chapno = -1 
\global\citationno=0
\global\headno = 0
\global\subheadno = 0
\global\procno = 0
\global\figno = 0

\def\resetcounters{%
\global\headno = 0%
\global\subheadno = 0%
\global\procno = 0%
\global\figno = 0%
}

\global\newchapflag=0 
\global\showpagenumflag=0 

\def\chinfo{\ifinappendices\alphanum\chapno\else\the\chapno\fi}%
\def\headinfo{\ifinappendices\alphanum\headno\else\the\headno\fi}%
\def\subheadinfo{\headinfo.\the\subheadno}
\def\procinfo{\headinfo.\the\procno}
\def\figinfo{\the\figno}        
\def\citationinfo{\the\citationno}%
\def\nextheadno{\global\advance\headno by 1 \global\subheadno = 0 \global\procno = 0}
\def\nextsubheadno{\global\advance\subheadno by 1}
\def\nextprocno{\global\advance\procno by 1 \procinfo}
\def\nextfigno{\global\advance\figno by 1 \figinfo}

{\global\let\noe=\noexpand%
%
%
\catcode`\@=11%
\catcode`\_=11%
\setcatcodes%
!global!def!_@@internal@@makeref#1{%
!global!expandafter!def!csname #1ref!endcsname##1{%
!csname _@#1@##1!endcsname%
!expandafter!ifx!csname _@#1@##1!endcsname!relax%
    !write16{#1 ##1 not defined - run saving references}%
    !undefinedreferencestrue%
!fi}}%
!global!def!_@@internal@@makelabel#1{%
!global!expandafter!def!csname #1label!endcsname##1{%
!edef!temptoken{!csname #1info!endcsname}%
!ifloadreferences%
    !expandafter!ifx!csname _@#1@##1!endcsname!relax%
        !write16{#1 ##1 not hitherto defined - rerun saving references}%
        !changedreferencestrue%
    !else%
        !expandafter!ifx!csname _@#1@##1!endcsname!temptoken%
        !else
            !write16{#1 ##1 reference has changed - rerun saving references}%
            !changedreferencestrue%
        !fi%
    !fi%
!else%
    !expandafter!edef!csname _@#1@##1!endcsname{!temptoken}%
    !edef!textoutput{!write!references{\global\def\_@#1@##1{!temptoken}}}%
    !textoutput%
!fi}}%
!global!def!makecounter#1{!_@@internal@@makelabel{#1}!_@@internal@@makeref{#1}}%
!unsetcatcodes%
}
\makecounter{ch}%
\makecounter{head}%
\makecounter{subhead}%
\makecounter{proc}%
\makecounter{fig}%
\makecounter{citation}%
\def\newref#1#2{%
\def\temptext{#2}%
\edef\bibliotextoutput{\expandafter\gobbleeight\meaning\temptext}%
\global\advance\citationno by 1\citationlabel{#1}%
\ifmakebiblio%
    \edef\fileoutput{\write\biblio{\noindent\hbox to 0pt{\hss$[\the\citationno]$}\hskip 0.2em\bibliotextoutput\medskip}}%
    \fileoutput%
\fi}%
\def\cite#1{%
$[\citationref{#1}]$%
\ifmakebiblio%
    \edef\fileoutput{\write\biblio{#1}}%
    \fileoutput%
\fi%
}%
%
%
%

\let\mypar=\par


\def\raggedleft{\leftskip=0pt plus 1fil \parfillskip=0pt}


\font\lettrinefont=cmr10 at 28pt
\def\lettrine #1[#2][#3]#4%
{\hangafter -#1 \hangindent #2
\noindent\hskip -#2 \vtop to 0pt{
\kern #3 \hbox to #2 {\lettrinefont #4\hss}\vss}}

\font\mylettrinefont=cmr10 at 28pt
\def\mylettrine #1[#2][#3][#4]#5%
{\hangafter -#1 \hangindent #2
\noindent\hskip -#2 \vtop to 0pt{
\kern #3 \hbox to #2 {\mylettrinefont #5\hss}\vss}}


\edef\Pagetitle={Blank}

\headline={\hfil\Pagetitle\hfil}

\footline={\hfil\myfontdefault\folio\hfil}

\def\nextoddpage
{
\newpage%
\ifodd\pageno%
\else%
    \global\showpagenumflag = 0%
    \null%
    \vfil%
    \eject%
    \global\showpagenumflag = 1%
\fi%
}


\def\newchap#1#2%
{%
%
%
\global\advance\chapno by 1%
\resetcounters%
%
%
\newpage%
\ifodd\pageno%
\else%
    \global\showpagenumflag = 0%
    \null%
    \vfil%
    \eject%
    \global\showpagenumflag = 1%
\fi%
\global\newchapflag = 1%
\global\showpagenumflag = 1%
%
%
{\font\chapfontA=cmsl10 at 30pt%
\font\chapfontB=cmsl10 at 25pt%
\null\vskip 5cm%
{\chapfontA\raggedleft\hfil%
{%
\ifnum\chapno=0
    \phantom{%
    \ifinappendices%
        Annexe \alphanum\chapno%
    \else%
        \the\chapno%
    \fi}%
\else%
    \ifinappendices%
        Annexe \alphanum\chapno%
    \else%
        \the\chapno%
    \fi%
\fi%
}%
\par}%
\vskip 2cm%
{\chapfontB\raggedleft%
\lineskiplimit=0pt%
\lineskip=0.8ex%
\hfil #1\par}%
\vskip 2cm%
}%
\edef\Pagetitle{#2}%
%
%
\ifmaketdm%
    \def\temp{#2}%
    \def\tempbis{\nobreak}%
    \edef\chaptitle{\expandafter\gobbleeight\meaning\temp}%
    \edef\mynobreak{\expandafter\gobbleeight\meaning\tempbis}%
    \edef\textoutput{\write\tdm{\bigskip{\noexpand\mytdmchapfont\noindent\chinfo\ - \chaptitle\hfill\noexpand\folio}\par\mynobreak}}%
\fi%
\textoutput%
}


\def\newhead#1%
{%
\ifhmode%
    \mypar%
\fi%
\ifnum\headno=0%
\ifinappendices
    \nobreak\vskip -\lastskip%
    \nobreak\vskip .5cm%
\fi
\else%
    \nobreak\vskip -\lastskip%
    \nobreak\vskip .5cm%
\fi%
\nextheadno%
\ifmaketdm%
    \def\temp{#1}%
    \edef\sectiontitle{\expandafter\gobbleeight\meaning\temp}%
    \edef\textoutput{\write\tdm{\noindent{\noexpand\mytdmheadfont\quad\headinfo\ - \sectiontitle\hfill\noexpand\folio}\par}}%
    \textoutput%
\fi%
\font\headfontA=cmbx10 at 14pt%
{\headfontA\noindent\headinfo\ - #1.\hfil}%
\nobreak\vskip .5cm%
}%


\def\newsubhead#1%
{%
\ifhmode%
    \mypar%
\fi%
\ifnum\subheadno=0%
\else%
    \penalty\headpenalty\vskip .4cm%
\fi%
\nextsubheadno%
\ifmaketdm%
    \def\temp{#1}%
    \edef\subsectiontitle{\expandafter\gobbleeight\meaning\temp}%
    \edef\textoutput{\write\tdm{\noindent{\noexpand\mytdmsubheadfont\quad\quad\subheadinfo\ - \subsectiontitle\hfill\noexpand\folio}\par}}%
    \textoutput%
\fi%
\font\subheadfontA=cmsl10 at 12pt
{\subheadfontA\noindent\subheadinfo\ #1.\hfil}%
\nobreak\vskip .25cm %
}%

%
%


\font\mathromanten=cmr10
\font\mathromanseven=cmr7
\font\mathromanfive=cmr5
\newfam\mathromanfam
\textfont\mathromanfam=\mathromanten
\scriptfont\mathromanfam=\mathromanseven
\scriptscriptfont\mathromanfam=\mathromanfive
\def\mathroman{\fam\mathromanfam}


\font\sf=cmss12

\font\sansseriften=cmss10
\font\sansserifseven=cmss7
\font\sansseriffive=cmss5
\newfam\sansseriffam
\textfont\sansseriffam=\sansseriften
\scriptfont\sansseriffam=\sansserifseven
\scriptscriptfont\sansseriffam=\sansseriffive
\def\mathsf{\fam\sansseriffam}


\font\bftwelve=cmb12

\font\boldten=cmb10
\font\boldseven=cmb7
\font\boldfive=cmb5
\newfam\mathboldfam
\textfont\mathboldfam=\boldten
\scriptfont\mathboldfam=\boldseven
\scriptscriptfont\mathboldfam=\boldfive
\def\mathbf{\fam\mathboldfam}


\font\mycmmiten=cmmi10
\font\mycmmiseven=cmmi7
\font\mycmmifive=cmmi5
\newfam\mycmmifam
\textfont\mycmmifam=\mycmmiten
\scriptfont\mycmmifam=\mycmmiseven
\scriptscriptfont\mycmmifam=\mycmmifive

\def\hexa#1{\ifcase #1 0\or 1\or 2\or 3\or 4\or 5\or 6\or 7\or 8\or 9\or A\or B\or C\or D\or E\or F\fi}
\mathchardef\mathi="7\hexa\mycmmifam7B
\mathchardef\mathj="7\hexa\mycmmifam7C


\font\mymsbmten=msbm10 at 8pt
\font\mymsbmseven=msbm7 at 5.6pt
\font\mymsbmfive=msbm5 at 4pt
\newfam\mymsbmfam
\textfont\mymsbmfam=\mymsbmten
\scriptfont\mymsbmfam=\mymsbmseven
\scriptscriptfont\mymsbmfam=\mymsbmfive

\mathchardef\mybeth="7\hexa\mymsbmfam69
\mathchardef\mygimmel="7\hexa\mymsbmfam6A
\mathchardef\mydaleth="7\hexa\mymsbmfam6B


\def\placelabel[#1][#2]#3{{%
\setbox10=\hbox{\raise #2cm \hbox{\hskip #1cm #3}}%
\ht10=0pt%
\dp10=0pt%
\wd10=0pt%
\box10}}%


\newif\ifinproclaim%
\global\inproclaimfalse%
\def\proclaim#1{%
\medskip%
%
%
\bgroup%
\inproclaimtrue%
\setbox10=\vbox\bgroup\leftskip=0.8em\noindent{\bftwelve #1}\sf%
}

\def\endproclaim{%
\egroup%
\setbox11=\vtop{\noindent\vrule height \ht10 depth \dp10 width 0.1em}%
\wd11=0pt%
\setbox12=\hbox{\copy11\kern 0.3em\copy11\kern 0.3em}%
\wd12=0pt%
\setbox13=\hbox{\noindent\box12\box10}%
\noindent\unhbox13%
\egroup%
\medskip\ignorespaces%
}

\def\proclaim#1{%
\medskip%
\bgroup%
\inproclaimtrue%
\noindent{\bftwelve #1}%
\nobreak\medskip%
\sf%
}

\def\endproclaim{%
\mypar\egroup\penalty\proclaimpenalty\medskip\ignorespaces%
}

\def\noskipproclaim#1{%
\medskip%
\bgroup%
\inproclaimtrue%
\noindent{\bf #1}\nobreak\sl%
}

\def\endnoskipproclaim{%
\mypar\egroup\penalty\proclaimpenalty\medskip\ignorespaces%
}


\def\ninn{{n\in\Bbb{N}}}

\def\proof{{\noindent\bf Proof:\ }}

\def\remark{{\noindent\sl Remark:\ }}

\def\mlim{\mathop{{\mathroman Lim}}}

\def\msf#1{{\mathsf #1}}

\def\qed{~$\square$}
\def\munion{\mathop{\cup}}
\def\minter{\mathop{\cap}}
\def\myitem#1{%
    \noindent\hbox to .5cm{\hfill#1\hss}
}

\catcode`\@=11
\def\Eqalign#1{\null\,\vcenter{\openup\jot\m@th\ialign{%
\strut\hfil$\displaystyle{##}$&$\displaystyle{{}##}$\hfil%
&&\quad\strut\hfil$\displaystyle{##}$&$\displaystyle{{}##}$%
\hfil\crcr #1\crcr}}\,}
\catcode`\@=12

\def\makeop#1{%
\global\expandafter\def\csname op#1\endcsname{{\mathroman #1}}}%

\def\makeopsmall#1{%
\global\expandafter\def\csname op#1\endcsname{{\mathroman{\lowercase{#1}}}}}%

\makeop{loc}
\makeop{Exp}
\makeop{Hess}
\makeopsmall{Tanh}%
\makeop{Id}

\font\mycirclefont=cmsy7
\def\textcircle{{\raise 0.3ex \hbox{\mycirclefont\char'015}}}

\let\emph=\bf

\hyphenation{quasi-con-formal}

%
%

\ifmakebiblio%
    \openout\biblio=biblio.tex %
    {%
        \edef\fileoutput{\write\biblio{\bgroup\leftskip=2em}}%
        \fileoutput
    }%
\fi%

\newref{BallSchGrom}{Ballmann W., Gromov M., Schroeder V., {\sl Manifolds of nonpositive curvature}, Progress in Mathematics, 61, Birkh\"auser Boston Inc., Boston, MA, (1985)}\newref{Brezis}{Brezis H., {\sl Functional analysis, Sobolev spaces and partial differential equations}, Universitext, Springer, New York, (2011)}
\newref{CaffNirSpruckV}{Caffarelli L., Nirenberg L., Spruck J., Nonlinear second-order elliptic equations. V. The Dirichlet problem for Weingarten hypersurfaces, {\sl Comm. Pure Appl. Math.} {\bf 41} (1988), no. 1, 47--70}
\newref{Gauss}{Gauss C. F., Disquisitiones generales circa superficies curvas, Oct. 8 1827}
\newref{Gromov}{Gromov M., Pseudoholomorphic curves in symplectic manifolds, {\sl Invent. Math.} {\bf 82} (1985), no. 2, 307--347}
\newref{GuanSpruck}{Guan B., Spruck J., The existence of hypersurfaces of constant Gauss curvature with prescribed boundary, {\sl J. Differential Geom.} {\bf 62} (2002), no. 2, 259--287}
\newref{ImaTani}{Imayoshi Y., Taniguchi M., {\sl An introduction to Teichm\"uller space}, Springer-Verlag, Tokyo, (1992)}
\newref{LabourieB}{Labourie F., Probl\`eme de Minkowski et surfaces \`a courbure constante dans les vari\'et\'es hyperboliques, {\sl Bull. Soc. Math. France} {\bf 119} (1991), no. 3, 307--325}
\newref{Labourie}{Labourie F., Un lemma de Morse pour les surfaces convexes, {\sl Invent. Math.} {\bf 141} (2000), no. 2, 239--297}
\newref{LehtoVirtanen}{Lehto O., Virtanen K. I., {\sl Quasiconformal mappings in the plane}, die Grundlehren der mathematischen Wissenschaften, {\bf 126}, Springer-Verlag, New-York, Heidelberg, (1973)}
\newref{Martinez}{Corro A. V., Mart\'inez A., Mil\'an F., Complete flat surfaces with two isolated singularities in hyperbolic 3-space, {\sl J. Math. Anal. Appl.} {\bf 366} (2010), no. 2, 582--592}
\newref{Pogorelov}{Pogorelov A. V., {\sl Extrinsic geometry of convex surfaces}, Translations of Mathematical Monographs, Vol. 35. American Mathematical Society, Providence, R.I., (1973)}
\newref{RosSpruck}{Rosenberg H., Spruck J., On the existence of convex hypersurfaces of constant Gauss curvature in hyperbolic space, {\sl J. Differential Geom.} {\bf 40} (1994), no. 2, 379--409}
\newref{Sasaki}{Sasaki S., On the differential geometry of tangent bundles of Riemannian manifolds, {\sl T\^ohoku Math. J. (2)} {\bf 10} (1958), 338--354}
\newref{SmiPKS}{Smith G., Pointed k-surfaces, {\sl Bull. Soc. Math. France} {\bf 134}, no. 4, (2006), 509--557}
\newref{SmiFCS}{Smith G., Moduli of Flat Conformal Structures of Hyperbolic Type, to appear in {\sl Geom. Dedicata}}
\newref{SmiPPH}{Smith G., Compactness for immersions of prescribed Gaussian curvature II - geometric aspects, arXiv:1002.2982}
\newref{TrudWang}{Trudinger N. S., Wang  X. J., On locally convex hypersurfaces with boundary, {\sl J. Reine Angew. Math.} {\bf 551} (2002), 11--32}

\ifmakebiblio%
    {\edef\fileoutput{\write\biblio{\egroup}}%
    \fileoutput}%
\fi%

%
%
%
\document
\myfontdefault
\global\chapno=1
\global\showpagenumflag=1
\def\Pagetitle{}
\null
\vfill
\def\centre{\rightskip=0pt plus 1fil \leftskip=0pt plus 1fil \spaceskip=.3333em \xspaceskip=.5em \parfillskip=0em \parindent=0em}%
\def\textmonth#1{\ifcase#1\or January\or Febuary\or March\or April\or May\or June\or July\or August\or September\or October\or November\or December\fi}
\font\abstracttitlefont=cmr10 at 14pt
{\abstracttitlefont\centre Hyperbolic Plateau problems\par}
\bigskip
{\centre Graham Smith\par}
\bigskip
{\centre 23rd May 2011\par}
\bigskip
{\centre Centre de Recerca Matem\`atica,\par
Facultat de Ci\`encies, Edifici C,\par
Universitat Aut\`onoma de Barcelona,\par
08193 Bellaterra,\par
Barcelona,\par
SPAIN\par}
\bigskip
\noindent{\emph Abstract:\ }We consider surfaces of constant Gaussian curvature immersed in $3$-dimensional manifolds, and we strengthen the compactness result of Labourie in the case where the ambient manifold is $3$-dimensional hyperbolic space. This allows us to prove results of existence of solutions to the asymptotic Plateau problem, as defined by Labourie, and the continuous dependence of these solutions on the data.
\par
\vfill
\nextoddpage
\def\Pagetitle{\sl Hyperbolic Plateau problems}
\global\pageno=1
\newhead{Introduction}
\noindent The Gaussian curvature of an immersed surface has interested differential geometers right from the very beginning with the announcement of Gauss' famous {\sl Teorema Egregium}, \cite{Gauss}. Having established this concept, it is natural to ask when there exist surfaces of constant Gaussian curvature in a given manifold subject to certain prescribed conditions which are typically geometric or topological in nature. For example, in the compact case, the Plateau problem and its various variants ask for compact immersed surfaces of constant Gaussian curvature with prescribed boundary (c.f. \cite{CaffNirSpruckV}, \cite{GuanSpruck}, \cite{RosSpruck}, \cite{TrudWang} etc.). In the non-compact case, Minkowski type problems ask for complete immersed hypersurfaces of constant Gaussian curvature whose normal vector field is prescribed by some sort of Gauss map (c.f. \cite{LabourieB}, \cite{Labourie}, \cite{Pogorelov} etc.).
\medskip
\noindent On closer examination, PDE considerations show that the qualitative nature of the problem is foremost affected by the sign of the Gaussian curvature. Indeed, the resulting PDE is elliptic, hyperbolic, or totally degenerate depending on whether the Gaussian curvature is positive, negative or zero respectively. Consequently, in the first case, the highly developed theory of regularity and compactness results for elliptic PDEs translates into a stronger geometric theory than can in general be obtained in the remaining cases. In the current setting, we go even further, as Labourie shows in \cite{Labourie} a straightforward but remarkable relationship between surfaces of constant positive Gaussian curvature on the one hand, and pseudo-holomorphic curves in the tangent bundle of the ambient manifold on the other which, through the powerful techniques developed by Gromov in \cite{Gromov}, readily yields a strong compactness result (Theorem \procref{ThmLabCompactness}) requiring weak hypotheses and involving a unique, elementary mode of degeneration which itself may often be excluded through straightforward geometric considerations. Using the continuity method and various limiting processes, this in turn yields existence results for solutions to the Plateau and Minkowksi problems (c.f. \cite{Labourie} and \cite{SmiPPH}).
\medskip
\noindent By definition, any compactness result assumes a choice of topology, and in the case of Labourie's compactness theorem (Theorem \procref{ThmLabCompactness}), the topology in question is the $C^\infty$ pointed Cheeger/Gromov topology on the space of complete, immersed surfaces (i.e. the topology of smooth convergence modulo reparametrisation). This topology is too weak to provide much control over the structure of the limit submanifolds. For example, the limit depends on the preferred points chosen, and in the case of surfaces, although on the one hand the genus is lower semi-continuous, the number of singularities is not even semi-continuous, compactness is not a closed property, and so on. As a consequence, it is rare that convergence in this topology is sufficient for any application, and supplementary techniques are often necessary in order to actually prove results. It is for this reason that in this paper, restricting attention to the case where the ambient manifold is $3$-dimensional hyperbolic space, we deduce a much stronger compactness theorem, which in turn permits us to recover a special case of the existence result of \cite{SmiFCS} as well as to show the continuous dependence of the solutions, which are unique, upon the data.
\medskip
\noindent Before stating our results, it is worth observing how the study of complete hypersurfaces of positive Gaussian curvature in $\Bbb{H}^3$ further subdivides into three different cases, depending on the value, $k$, of the curvature. When $k>1$, the intrinsic geometry is spherical, and we readily show that the immersion is a geodesic sphere of radius determined by $k$. When $k=1$, the intrinsic geometry of the surface is Euclidean, and the immersion is either a horosphere, or a cylinder of points lying at constant distance from a complete geodesic. It is thus only in the case where $k\in]0,1[$, and the intrinsic geometry is hyperbolic, that a rich zoology arises, and it is for this reason that we restrict attention henceforth to these values of $k$ (note, however, that by admitting singularities nonetheless we recover an interesting zoology also in the preceding two cases, as in for example \cite{Martinez}).
\medskip
\noindent We now describe the setting in which our results are proven. Let $\Bbb{H}^3$ be $3$-dimensional hyperbolic space. Let $T\Bbb{H}^3$ be the tangent bundle over $\Bbb{H}^3$ and let $U\Bbb{H}^3$ be the circle bundle of vectors of unit length in $T\Bbb{H}^3$. Let $S$ be an oriented surface, let $i:S\rightarrow\Bbb{H}^3$ be an immersion, and let $\msf{N}_i:S\rightarrow\Bbb{H}^3$ be the unit normal vector field over $i$ compatible with the orientation. We define the {\bf Gauss lift} $\hat{\mathi}:S\rightarrow U\Bbb{H}^3$ by:
$$
\hat{\mathi} = \msf{N}_i.
$$
\proclaim{Definition \nextprocno}
\noindent For $k\in]0,1[$, we say that the immersed hypersurface $\Sigma:=(S,i)$ is a $k$-surface if and only if:
\medskip
\myitem{(i)} $i$ is everywhere locally strictly convex and is of constant Gaussian curvature equal to $k$; and
\medskip
\myitem{(ii)} $\hat{\mathi}$ is complete.
\endproclaim
\remark The condition that the Gauss lift rather than the immersion itself be complete, although apparently strange, is natural from the perspective of Labourie's compactness theorem which, strictly speaking, operates entirely within the unitary bundle (c.f. \cite{Labourie}). This often leads to subtle complexities in the geometric arguments required. However, in the case of most interest to us, being that of solutions generated by ramified coverings of the sphere, completeness of $i$ follows from completeness of $\hat{\mathi}$ (c.f. \cite{SmiPKS}).
\medskip
\noindent Let $\partial_\infty\Bbb{H}^3$ be the ideal boundary of $\Bbb{H}^3$ (c.f. \cite{BallSchGrom}). We define the {\bf Gauss-Minkowksi} map, $\overrightarrow{n}:U\Bbb{H}^3\rightarrow\partial_\infty\Bbb{H}^3$ as follows: choose $V\in U\Bbb{H}^3$, and define $\overrightarrow{n}(V)$ to be the point in $\partial_\infty\Bbb{H}^3$ towards which $V$ points. Explicitly, if $\gamma:\Bbb{R}\rightarrow\Bbb{H}^3$ is the unique geodesic such that:
$$
(\partial_t\gamma)(0) = V,
$$
\noindent then:
$$
\overrightarrow{n}(V) = \gamma(+\infty) = \mlim_{t\rightarrow+\infty}\gamma(t).
$$
\noindent We recall that $\partial_\infty\Bbb{H}^3$ has the conformal structure of the Riemann sphere. Following \cite{Labourie}, we make the following definition:
\goodbreak
\proclaim{Definition \nextprocno}
\noindent An asymptotic Plateau problem is a pair $(S,\varphi)$ where $S$ is a Riemann surface, and $\varphi:S\rightarrow\Bbb{H}^3$ is a locally conformal mapping.
\medskip
\noindent Given an asymptotic Plateau problem, $(S,\varphi)$, and given $k\in]0,1[$, the smooth immersion $i:S\rightarrow\Bbb{H}^3$ is said to be a solution to $(S,\varphi)$ if and only if:
\medskip
\myitem{(i)} $(S,i)$ is a $k$-surface; and
\medskip
\myitem{(ii)} $\varphi = \overrightarrow{n}\circ i$.
\endproclaim
\remark This definition superficially appears more restrictive than that given in \cite{Labourie}, where $S$ may be any topological surface and $\varphi$ a local homeomorphism. However, pulling back the conformal structure of $\partial_\infty\Bbb{H}^3$ through $\varphi$ readily transforms an asymptotic Plateau problem in the sense of \cite{Labourie} into an asymptotic Plateau problem in the current sense, and so we see that no generality is lost.
\medskip
\remark we see that an asymptotic Plateau problem is in fact closer to the Minkowski problem than the Plateau problem. We nonetheless retain the current terminology for consistency with \cite{Labourie} and \cite{SmiPKS}.
\medskip
\noindent The main result of this paper is the following Theorem which follows immediately from Proposition \procref{PropStrongCompactness}:
\proclaim{Theorem \nextprocno, {\bf Compactness}}
\noindent Let $\Bbb{D}$ be the Poincar\'e disc. Let $(\varphi_n)_\ninn:\Bbb{D}\rightarrow\partial_\infty\Bbb{H}^3$ be locally conformal maps. Choose $k\in]0,1[$ and for all $n$ let $i_n:\Bbb{D}\rightarrow\Bbb{H}^3$ be a $k$-surface which is a solution to the asymptotic Plateau problem $(\Bbb{D},\varphi_n)$. Suppose that there exists a locally conformal map $\varphi_0:\Bbb{D}\rightarrow\partial_\infty\Bbb{H}^3$ towards which $(\varphi_n)_\ninn$ converges in the local uniform sense. Then there exists an immersion $i_0:\Bbb{D}\rightarrow\Bbb{H}^3$ such that
\medskip
\myitem{(i)} $(\Bbb{D},i_0)$ is a $k$-surface;
\medskip
\myitem{(ii)} $i_0$ is a solution to the asymptotic Plateau problem $(\Bbb{D},\varphi_0)$; and
\medskip
\myitem{(ii)} $(i_n)_\ninn$ converges to $i_0$ in the $C^\infty_\oploc$ sense.
\endproclaim
\proclabel{ThmCompactness}
\noindent This yields a short, alternative proof of the following special case of Theorem $1.1$ of \cite{SmiFCS}:
\proclaim{Theorem \nextprocno, {\bf Existence}}
\noindent Let $\Bbb{D}$ be the Poincar\'e disc. Let $\varphi:\Bbb{D}\rightarrow\partial_\infty\Bbb{H}^3$ be a locally conformal map. For all $k\in]0,1[$, there exists a unique immersion $i:\Bbb{D}\rightarrow\Bbb{H}^3$ such that:
\medskip
\myitem{(i)} $(\Bbb{D},i)$ is a $k$-surface; and
\medskip
\myitem{(ii)} $i$ is a solution to the asymptotic Plateau problem $(\Bbb{D},\varphi)$.
\endproclaim
\proclabel{ThmExistence}
\noindent Moreover the following continuity result, being valid also in the non-cocompact case, is thus stronger than what may be obtained via the techniques of \cite{SmiFCS}:
\proclaim{Theorem \nextprocno, {\bf Continuous Dependence}}
\noindent Let $\Bbb{D}$ be the Poincar\'e disc. Let $(\varphi_n)_\ninn, \varphi_0:\Bbb{D}\rightarrow\partial_\infty\Bbb{H}^3$ be locally conformal maps. Choose $k\in]0,1[$, and for all $n\in\Bbb{N}\munion\left\{0\right\}$, let $i_n:\Bbb{D}\rightarrow\Bbb{H}^3$ be the unique immersion such that:
\medskip
\myitem{(i)} $(\Bbb{D},i_n)$ is a $k$-surface; and
\medskip
\myitem{(ii)} $i_n$ is the unique solution to the asymptotic Plateau problem $(\Bbb{D},\varphi)$.
\medskip
\noindent Suppose that $(\varphi_n)_\ninn$ converges to $\varphi_0$ in the local uniform sense, then $(i_n)_\ninn$ converges to $i_0$ in the $C^\infty_\oploc$ sense.
\endproclaim
\proclabel{ThmContinuousDependance}
\noindent This paper is a greatly revised version of the second chapter of the author's doctoral thesis. The author would like to thank Fran\c{c}ois Labourie for having proposed this problem, for his guidance during that period, and for his encouragement to prepare the current version. The author would also like to thank the Universit\'e Paris XI, and the Max Planck Insitute for Mathematics in the Sciences in Leipzig for providing the conditions required to prepare the previous version of this paper. The author would like thank the Centre de Recerca Matem\`atica in Barcelona for providing the conditions required to prepare the current version of this paper which was written whilst the author was benifitting from a Marie Curie Postdoctoral fellowship.
\newhead{Conformal Structure}
\noindent For $r>0$, let $S_r\Bbb{H}^3$ be the bundle of spheres of radius $r$ in $T\Bbb{H}^3$. Let $S$ be an oriented surface and let $i:S\rightarrow\Bbb{H}^3$ be an immersion. Let $\msf{N}_i$ be the unit normal vector field over $i$ compatible with the orientation. We define $\hat{\mathi}_r:S\rightarrow S_r\Bbb{H}^3$ by:
$$
\hat{\mathi}_r = r\msf{N}_i.
$$
\noindent We call $\hat{\mathi}_r$ the {\bf Gauss lift} of $i$ at height $r$. Let $\hat{g}$ be the Sasaki metric over $T\Bbb{H}^3$ (c.f. \cite{Sasaki}). Let $I_i$ and $III_i$ be the first and third fundamental forms of $i$. We readily see that:
$$
\hat{\mathi}_r^*\hat{g} = I_i + r^2III_i.
$$
\noindent Let $J_{i,r}$ be the complex structure compatible with the orientation induced over $\Sigma$ by the metric $\hat{\mathi}_r^*\hat{g}$. We henceforth refer to $J_{i,r}$ as the complex structure {\bf induced} by $\hat{\mathi}_r$. We shall see that this complex structure plays a central role in the sequel.
\medskip
\noindent Let $\partial_\infty\Bbb{H}^3$ be the ideal boundary of $\Bbb{H}^3$. We define the {\bf Gauss-Minkowski} mapping $\overrightarrow{n}:S_r\Bbb{H}^3\rightarrow\partial_\infty\Bbb{H}^3$ as in the introduction. Given an immersion $i:S\rightarrow\Bbb{H}^3$, the composition $(\overrightarrow{n}\circ\hat{\mathi}_{i,r})$ (which is indepedant of $r$) is everywhere an orientation preserving local diffeomorphism between the Riemann surfaces $(S,J_{i,r})$ and $\partial_\infty\Bbb{H}^3=\hat{\Bbb{C}}$. In general, this map is not conformal but is quasiconformal. We thus recall the definition of quasiconformality. Let $U\subseteq\Bbb{C}$ be open and let $\alpha:U\rightarrow\Bbb{C}$ be a smooth, orientation preserving mapping which is everywhere a local diffeomorphism. For all $z\in U$, the derivative of $\alpha$ is given by:
$$
D\alpha_z\cdot w = \partial\alpha(z)w + \overline{\partial}\alpha(z)\overline{w}.
$$
\noindent We define $\mu(\alpha)$, the {\bf complex dilatation} of $\alpha$, by:
$$
\mu(\alpha)(z) = \frac{\overline{\partial}\alpha(z)}{\partial\alpha(z)}
$$
\noindent For $K\geqslant 1$, we say that $\alpha$ is $K$-{\bf quasiconformal} if and only if for all $z\in U$:
$$
\left|\mu(\alpha)(z)\right| \leqslant \frac{K-1}{K+1}.
$$
\noindent We recall that the concept of $K$-quasiconformality readily extends to orientation preserving local homeomorphisms, and that a local homeomorphism is $1$-quasiconformal if and only if it is smooth and conformal (c.f. \cite{ImaTani}).
\proclaim{Lemma \nextprocno}
\noindent Choose $k\in]0,1[$ and let $r=k^{-1/2}>1$. Let $i:\Sigma\rightarrow\Bbb{H}^3$ be a locally strictly convex immersion, let $\hat{\mathi}_r:S\rightarrow S_r\Bbb{H}^3$ be the Gauss lift of $i$ at height $r$, and let $J_{i,r}$ be the complex structure induced over $S$ by $\hat{\mathi}_r$. If $i$ has constant Gaussian curvature equal to $k$, then $(\overrightarrow{n}\circ\hat{\mathi}_r):(S,J_{i,r})\rightarrow\partial_\infty\Bbb{H}^3$ is $r$-quasiconformal.
\endproclaim
\proclabel{LemmaQuasiConformal}
\proof We identify $\Bbb{H}^3$ with $H^+$, the upper half space in $\Bbb{R}^3$. We identify $S_r\Bbb{H}^3$ with $H^+\times S^2$. Using elementary Euclidean geometry, we readily show:
$$
\overrightarrow{n}((x,y,z),(u,v,w)) = (x,y) + \frac{z}{1-w}(u,v).
$$
\noindent Choose $p_0\in S$. By applying an isometry of $\Bbb{H}^3$, we may assume that there exists a neighbourhood, $U$, of $0$ in $\Bbb{R}^2$, a neighbourhood, $V$, of $p_0$ in $S$ and a smooth function $f:\Bbb{R}^2\rightarrow]0,\infty[$ such that $i(V)$ coincides with the graph of $f$ over $U$. We may suppose, moreover, that $f(0)=1$ and $Df(0)=0$. Since the result is invariant under reparametrisation, it suffices to prove it for $\hat{f}:U\rightarrow\Bbb{H}^3$ given by:
$$
\hat{f}(x,y) = (x,y,f(x,y)).
$$
\noindent Let $(e_1,e_2)$ be an orthonormal basis of $\Bbb{R}^2$ with respect to which $\opHess_0(f)$ is diagonal, and let $(\lambda_1,\lambda_2)$ be the corresponding eigenvalues. Let $\msf{N}_{\hat{f}}$ be the unit normal vector field over $\hat{f}$ compatible with the orientation. We may assume that $\msf{N}_{\hat{f}}$ points downwards. Since the hyperbolic metric over $H^+$ is conformally equivalent to the Euclidean metric, $\msf{N}_{\hat{f}}$ coincides up to a scaling factor with the Euclidean unit normal, and we readily obtain:
$$
\msf{N}_{\hat{f}}(x,y) = ((x,y,f(x,y)),(1+\|D f\|^2)^{-1}(\partial_xf,\partial_yf,-1))
$$
\noindent  Bearing in mind that $\hat{f}(0)=(0,0,0)$ and $Df(0)=(0,0)$, using the chain rule, we readily obtain:
$$
D(\overrightarrow{n}\circ\msf{N}_{\hat{f}})(0) = M := \opId + \frac{1}{2}\pmatrix \lambda_1\hfill& 0\hfill\cr 0\hfill&\lambda_2\hfill\cr\endpmatrix.
$$
\noindent Let $I_0(f)$, $II_0(f)$ and $III_0(f)$ be the first, second and third fundamental forms of $\hat{f}$ with respect to the hyperbolic metric at $0$. Since $f(0)=1$ and $Df(0)=0$, we readily see that $I$ coincides with the Euclidean metric over $\Bbb{R}^2$. Let $II_0^e(f)$ be the second fundamental form of $\hat{f}$ with respect to the Euclidean metric at $0$. Observe that $II_0(f)-II_0^e(f)$ only depends on $\msf{N}_{\hat{f}}(0)$. Let $g:U\rightarrow\Bbb]0,\infty[$ be the function whose graph is a portion of the unit sphere in $\Bbb{R}^3$ centred at $0$. In particular, $g(0)=1$ and $Dg(0)=0$. We choose $\msf{N}_{\hat{g}}(0)$ to be the downward pointing normal. Then, by definition:
$$
II_0^e(g)(e_i,e_j) = -\delta_{ij},
$$
\noindent However, the graph of $g$ is a totally geodesic subspace of $\Bbb{H}^3$. Thus:
$$
II_0(g) = 0.
$$
\noindent Since $\msf{N}_{\hat{f}}$ points downwards, we readily obtain:
$$
II_0^e(f)(e_i,e_j) = \opHess_0(f)(e_i,e_j).
$$
\noindent Thus:
$$\matrix
II_0(f)(e_i,e_j)\hfill&=II_0^e(f)(e_i,e_j) + (II_0(f) - II_0^e(f))(e_i,e_j)\hfill\cr
&=\partial_i\partial_jf(0) + \delta_{ij}.\hfill\cr
\endmatrix$$
\noindent In summary:
$$
I_0(f) = \pmatrix 1\hfill& 0\hfill\cr 0\hfill&1\hfill\cr\endpmatrix,\quad
II_0(f) = \pmatrix 1 + \lambda_1\hfill& 0\hfill\cr 0\hfill&1+\lambda_2\hfill\cr\endpmatrix,\quad
III_0(f) = \pmatrix (1 + \lambda_1)^2\hfill& 0\hfill\cr 0\hfill&(1+\lambda_2)^2\hfill\cr\endpmatrix.
$$
\noindent Using the formula for the Gaussian curvature of $\hat{f}$, we obtain:
$$
(1+\lambda_1)(1+\lambda_2) = k.
$$
\noindent Moreover, if $J$ is the conformal structure induced over $S$ by the Gauss lift at height $r$ of $\hat{f}$, then:
$$
Je_1 = r(1 + \lambda_1)e_2,\qquad Je_2 = -r(1+\lambda_2)e_1.
$$
\noindent We thus define $(e_1',e_2')$ by:
$$
e_1' = \frac{1}{\sqrt{1+\lambda_1}} e_1,\qquad e_2'=\frac{1}{\sqrt{1+\lambda_2}}e_2.
$$
\noindent So that $Je_1'=e_2'$ and $Je_2'=-e_1'$. With respect to the basis $(e_1',e_2')$ in the domain, and the basis $(e_1,e_2)$ in the range, the matrix of $D(\overrightarrow{n}\circ\msf{N}_{\hat{f}})$ is therefore equal to $M'$, where:
$$
M' = \frac{1}{2}\pmatrix \frac{2+\lambda_1}{\sqrt{1+\lambda_1}}\hfill& 0\hfill\cr 0\hfill&\frac{2+\lambda_2}{\sqrt{1+\lambda_2}}\hfill\cr\endpmatrix.
$$
\noindent The complex dilatation of $M'$ is equal to $\mu$, where:
$$
\mu = \frac{(\sqrt{1+\lambda_1} - \sqrt{1+\lambda_2})}{(\sqrt{1+\lambda_1} + \sqrt{1+\lambda_2})}\frac{r-1}{r+1}.
$$
\noindent The absolute value of the first multiplicand is trivially no greater than $1$, and so:
$$
\left|\mu\right| \leqslant \frac{r-1}{r+1}.
$$
\noindent $M'$ is therefore $r$-quasiconformal, and this completes the proof.\qed
\newhead{Weak Compactness}
\noindent We first discuss the form of degenerate limit that may arise. Let $\Gamma\subseteq\Bbb{H}^3$ be a complete geodesic in $\Bbb{H}^3$. Let $\msf{N}_r(\Gamma)\subseteq S_r\Bbb{H}^3$ be the set of vectors of length $r$ over $\Gamma$ which are normal to $\Gamma$.
\proclaim{Definition \nextprocno}
\noindent Let $S$ be a surface and let $\hat{\mathi}:S\rightarrow S_r\Bbb{H}^3$ be a complete immersion. We say that $(S,\hat{\mathi})$ is tubular if and only if there exists a complete geodesic $\Gamma$ such that $\hat{\mathi}$ is a covering map of $\msf{N}_r(\Gamma)$.
\endproclaim
\remark Observe that if $(S,\hat{\mathi})$ is tubular, and if $J$ is the complex structure induced over $S$ by $\hat{\mathi}$, then $(S,J)$ is of parabolic type. This plays an important role in the sequel.
\medskip
\noindent We recall the definition of convergence modulo reparametrisation for complete immersed submanifolds. Let $(M,g)$ be a Riemannian manifold. For all $n\in\Bbb{N}$, let $i_n:S_n\rightarrow M$ be a complete immersion, and for all $n$, choose $p_n\in S_n$. Let $i_0:S_0\rightarrow M$ be another complete immersion and choose $p_0\in S_0$.
\proclaim{Definition \nextprocno}
\noindent We say that $(S_n,i_n,p_n)_\ninn$ converges to $(S_0,i_0,p_0)$ in the $C^\infty$ sense modulo repara\-metrisation if and only if there exists a sequence of mappings $\varphi_n:S_0\rightarrow S_n$ such that:
\medskip
\myitem{(i)} for all $n$, $\varphi_n(p_0)=p_n$; and
\medskip
\noindent for every relatively compact open subset $\Omega\subseteq S_0$, there exists $N\in\Bbb{N}$ such that:
\medskip
\myitem{(ii)} for all $n\geqslant N$, the restriction of $\varphi_n$ to $\Omega$ is a diffeomorphism onto its image; and
\medskip
\myitem{(iii)} $(i_n\circ\varphi_n)_\ninn$ converges to $i_0$ in the $C^\infty_\oploc$ sense over $\Omega$.
\medskip
\noindent We will refer to $(\varphi_n)_\ninn$ as a sequence of convergence mappings for $(S_n,i_n,p_n)_\ninn$ with respect to $(S_0,i_0,p_0)$.
\endproclaim
\proclabel{DefnCheegerGromov}
\noindent Choose $k\in]0,1[$ and let $r=k^{-1/2}$. For $n\in\Bbb{N}$, let $S_n$ be an oriented surface, let $i_n:S_n\rightarrow\Bbb{H}^3$ be a locally strictly convex immersion of constant Gaussian curvature equal to $k$, and let $p_n$ be a point in $S_n$. For all $n$, let $\hat{\mathi}_{n,r}:S_n\rightarrow S_r\Bbb{H}^3$ be the Gauss lift of $i_n$ at height $r$. In \cite{Labourie}, Labourie proves:
\goodbreak
\proclaim{Theorem \nextprocno}
\noindent Suppose that:
\medskip
\myitem{(i)} for all $n$, $\hat{\mathi}_{n,r}:S\rightarrow S_r\Bbb{H}^3$ is a complete immersion; and
\medskip
\myitem{(ii)} there exists a compact subset $K\subseteq S_r\Bbb{H}^3$ such that $\hat{\mathi}_{n,r}(p_n)\in K$ for all $n$.
\medskip
\noindent Then there exists a complete pointed immersed surface $(S_0,\hat{\mathi}_0,p_0)$ towards which $(S_n,\hat{\mathi}_n,p_n)\subseteq S_r\Bbb{H}^3$ subconverges in the $C^\infty$ sense modulo reparametrisation. Moreover either:
\medskip
\myitem{(i)} $\pi\circ\hat{\mathi}_0$ is an immersion, where $\pi:S_r\Bbb{H}^3\rightarrow\Bbb{H}^3$ is the canonical projection; or
\medskip
\myitem{(ii)} $(S_0,\hat{\mathi}_0)$ is tubular.
\endproclaim
\proclabel{ThmLabCompactness}
\remark Observe that, in the first case, $\pi\circ\hat{\mathi}_0$ is an immersion of constant Gaussian curvature equal to $k$. Importantly, it is not necessarily complete.
\medskip
\noindent We now relate the convergence maps to the conformal structure of the limit. We require the following Schwartz lemma:
\proclaim{Lemma \nextprocno}
\noindent Let $S$ be an oriented surface. Let $(J_n)_\ninn$ be a sequence of smooth almost complex structures over $S$ which converges to another almost complex structure $J_0$ in the $C^\infty_\oploc$ sense. For all $n\in\Bbb{N}$, let $\alpha_n:S\rightarrow\Bbb{D}$ be holomorphic with respect to $J_n$. There exists $\alpha_0:S\rightarrow\Bbb{D}$ towards which $(\alpha_n)_\ninn$ subconverges in the $C^\infty_\oploc$ sense.
\endproclaim
\proclabel{LemmaSchwartz}
\proof Choose $p\in S$. Let $U$ be a neighbourhood of $p$, and for all $n\in\Bbb{N}\munion\left\{0\right\}$, let $(e_{1,n},e_{2,n})$ be a frame defined over $U$ such that:
\medskip
\myitem{(i)} for all $n$, $J_ne_{1,n}=e_{2,n}$; and
\medskip
\myitem{(ii)} $(e_{1,n},e_{2,n})_\ninn$ converges to $(e_{1,0},e_{2,0})$ in the $C^\infty_\oploc$ sense.
\medskip
\noindent For all $n\in\Bbb{N}\munion\left\{0\right\}$, define $a_n,b_n:\Omega\rightarrow\Bbb{R}$ such that:
$$
[e_{1,n},e_{2,n}] = a_ne_{1,n} + b_ne_{2,n}.
$$
\noindent We define the generalised Laplacian $L_n$ such that, for all $f$:
$$
L_nf = (e_{1,n}e_{1,n} + e_{2,n}e_{2,n})f - (b_ne_{1,n} - a_ne_{2,n})f.
$$
\noindent For all $n\in\Bbb{N}$, we denote $\alpha_n=\xi_n + i\eta_n$. The Cauchy-Riemann equations yield:
$$\matrix
(e_{1,n}e_{1,n}+e_{2,n}e_{2,n})\xi_n \hfill&= (e_{1,n}e_{2,n} - e_{2,n}e_{1,n})\eta_n\hfill\cr
&= (a_ne_{1,n} + b_ne_{2,n})\eta_n\hfill\cr
&= (b_ne_{1,n} - a_ne_{2,n})\xi_n.\hfill\cr
\endmatrix$$
\noindent Thus:
$$
L_n\xi_n=0.
$$
\noindent For all $n$:
$$
\|\xi_n\|_{L^\infty}\leqslant \|\alpha_n\|_{L^\infty}\leqslant 1.
$$
\noindent It thus follows by classical elliptic regularity (c.f. \cite{Brezis}) that there exists $\xi_0$ towards which $(\xi_n)_\ninn$ subconverges. Likewise, there exists $\eta_0$ towards which $(\eta_n)_\ninn$ subconverges. This completes the proof.\qed
\medskip
\noindent We refine Labourie's result as follows: choose $k\in]0,1[$ and denote $r=k^{-1/2}$. Let $\Bbb{D}$ be the Poincar\'e disc. For all $n\in\Bbb{N}$, let $i_n:\Bbb{D}\rightarrow\Bbb{H}^3$ be a locally strictly convex immersion of constant Gaussian curvature equal to $k$, let $\hat{\mathi}_{n,r}$ be its Gauss lift at height $r$. Suppose that $\hat{\mathi}_{n,r}:\Bbb{D}\rightarrow S_r\Bbb{H}^3$ is complete, and, moreover, for all $n$, the complex structure induced over $\Bbb{D}$ by $\hat{\mathi}_{n,r}$ coincides with the canonical complex structure of $\Bbb{D}$.
\medskip
\noindent For all $n\in\Bbb{N}$, choose $z_n\in\Bbb{D}$. Let $K\subseteq S_r\Bbb{H}^3$ be compact, and suppose that $\hat{\mathi}_{n,r}(z_n)\in K$. By Theorem \procref{ThmLabCompactness}, there exists a complete immersed surface $(S_0,\hat{\mathi}_0,p_0)$ in $S_r\Bbb{H}^3$ towards which $(\Bbb{D},\hat{\mathi}_{n,r},z_n)$ subconverges modulo reparametrisation in the $C^\infty_\oploc$ sense. Let $(\varphi_n)_\ninn:S_0\rightarrow\Bbb{D}$ be a sequence of convergence mappings of $(\Bbb{D},\hat{\mathi}_{n,r},z_n)$ with respect to $(S_0,\hat{\mathi}_0,p_0)$. Let $J_0$ be the complex structure induced over $S_0$ by $\hat{\mathi}_0$.
\proclaim{Lemma \nextprocno}
\noindent There exists a conformal mapping $\varphi_0:(S_0,J_0)\rightarrow\overline{\Bbb{D}}$ towards which $(\varphi_n)_\ninn$ subconverges in the $C^\infty_\oploc$ sense.
\endproclaim
\proclabel{LemmaConvergenceMappingsConverge}
\remark Recall the remark following Definition \procref{DefnCheegerGromov}. If $(S,J_0)$ is tubular, then, in particular it is of parabolic type, and so $\varphi_0$ is constant.
\medskip
\proof Let $\Omega\subseteq S_0$ be a relatively compact open set. Choose $N\in\Bbb{N}$ such that, for $n\geqslant N$, the restriction of $\varphi_n$ to $\Omega$ is a diffeomorphism onto its image. For all $n$, let $J_n$ be the complex structure induced over $S_0$ by $(\hat{\mathi}_{n,r}\circ\varphi_n)$. Since $(\hat{\mathi}_{n,r})_{n\geqslant N}$ subconverges to $\hat{\mathi}_0$ in the $C^\infty_\oploc$ sense over $\Omega$, $(J_n)_\ninn$ also subconverges to $J_0$ in the $C^\infty_\oploc$ sense over $\Omega$. However, for all $n$, by definition $\varphi_n:(\Omega,J_n)\rightarrow\Bbb{D}$ is conformal. It follows by Lemma \procref{LemmaSchwartz} that there exists a conformal map $\varphi_0:(\Omega,J_0)\rightarrow\Bbb{D}$ towards which $(\varphi_n)_\ninn$ subconverges in the $C^\infty_\oploc$ sense. The result now follows by a diagonal argument.\qed
\newhead{Compactness of Quasiconformal Mappings}
\noindent Let $\Bbb{D}$ be the Poincar\'e disc. We identify $\partial_\infty\Bbb{H}^3$, the ideal boundary of hyperbolic space, with the Riemann sphere, $\hat{\Bbb{C}}$. Let $(\varphi_n)_\ninn,\varphi_0:\Bbb{D}\rightarrow\hat{\Bbb{C}}$ be locally conformal maps such that $(\varphi_n)_\ninn$ converges to $\varphi_0$ in the $C^\infty_\oploc$ sense.
\medskip
\noindent Choose $k\in]0,1[$ and denote $r=k^{-1/2}$. For all $n$, let $i_n:\Bbb{D}\rightarrow\Bbb{H}^3$ be an immersion of constant Gaussian curvature equal to $k$ such that:
\medskip
\myitem{(i)} the Gauss lift, $\hat{\mathi}_{n,r}$, of $i$ at height $r$ is a complete immersion of $\Bbb{D}$ into $S_r\Bbb{H}^3$; and
\medskip
\myitem{(ii)} $(\overrightarrow{n}\circ\hat{\mathi}_{n,r})=\varphi_n$.
\medskip
\noindent For all $n$, let $J_n$ be the complex structure induced over $\Bbb{D}$ by $\hat{\mathi}_{n,r}$. Let $\mu_n$ be the complex dilatation of $J_n$ with respect to the canonical complex structure over $\Bbb{D}$. By Lemma \procref{LemmaQuasiConformal}, for all $n$:
$$
\|\mu_n\|\leqslant\frac{r-1}{r+1}.
$$
\proclaim{Proposition \nextprocno}
\noindent For all $n$, there exists a conformal homeomorphism $\alpha_n:(\Bbb{D},J_n)\rightarrow\Bbb{D}$ such that $\alpha_n(0)=0$.
\endproclaim
\proclabel{PropConformalRepresentation}
\proof Since $J_n$ is $r$-quasiconformal with respect to the canonical complex structure over $\Bbb{D}$, $(\Bbb{D},J_n)$ is of hyperbolic type (c.f. \cite{ImaTani}). The result follows by Riemann's uniformisation theorem.\qed
\proclaim{Proposition \nextprocno}
\noindent There exists a quasiconformal mapping $\alpha_0:\Bbb{D}\rightarrow\Bbb{D}$ such that $\alpha_0(0)=0$ and $(\alpha_n)_\ninn$ subconverges to $\alpha_0$ in the $C^0$ sense.
\endproclaim
\proclabel{PropQuasiConfCompactness}
\proof We apply a M\"obius transformation and work in the upper half space, $\Bbb{H}^+$. Choose $n\in\Bbb{N}$. Conjugating by a M\"obius transformation, we may suppose that $\alpha_n$ is an $r$-quasiconformal homeomorphism from $\Bbb{H}^+$ to itself which fixes $i$. We extend $\alpha_n$ to an $r$-quasiconformal homeomorphism of $\hat{\Bbb{C}}$ as follows: let $\mu_n$ be the complex dilatation of $\alpha_n$. We extend $\mu_n$ uniquely to a bounded, measurable $(1,1)$-form over $\Bbb{C}$ such that, for all $z\in\Bbb{H}^+$:
$$
\mu_n(\overline{z}) = \overline{\mu}_n(z).
$$
\noindent By Theorem $4.25$ of \cite{ImaTani}, there exists a unique quasiconformal homeomorphism $\tilde{\alpha}_n:\hat{\Bbb{C}}\rightarrow\hat{\Bbb{C}}$ preserving $0$, $1$ and $+\infty$ and whose conformal dilatation is equal to $\mu_n$. Since its conformal dilatation and fixed points are symmetric about the real axis, $\tilde{\alpha}_n$ preserves $\Bbb{R}$. Since it fixes $0$, $1$ and $\infty$, it also preserves the orientation of $\Bbb{R}$, and so, since it is orientation preserving, it maps $\Bbb{H}^+$ to itself.
\medskip
\noindent For all $n$, let $\beta_n:\hat{\Bbb{C}}\rightarrow\hat{\Bbb{C}}$ be a M\"obius mapping also preserving $\Bbb{H}^+$ such that:
$$
\beta_n(\tilde{\alpha}_n(i)) = i.
$$
\noindent Denote $\alpha_n'=\beta_n\circ\tilde{\alpha}_n$. The mapping $\alpha_n$ is a homeomorphism of $\Bbb{H}^+$ preserving $i$. Moreover, its conformal dilatation over $\Bbb{H}^+$ is equal to that of $\tilde{\alpha}_n$, which is in turn equal to that of $\alpha_n$. It follows that, after composing $\beta_n$ with a rotation if necessary, $\alpha_n'=\alpha_n$ over $\Bbb{H}^+$ (c.f. \cite{ImaTani}), and so $\alpha_n'$ is the desired extension of $\alpha_n$. Henceforth, we denote $\alpha_n'$ merely by $\alpha_n$. Observe that, for all $z\in\Bbb{H}^+$:
$$
\alpha_n(\overline{z}) = \overline{\alpha}_n(z).
$$
\noindent In particular, $\alpha_n$ preserves $i$, $-i$ and $\Bbb{R}$. Thus, by Theorem $II.5.1$ of \cite{LehtoVirtanen}, the sequence $(\alpha_n)_\ninn$ constitutes a normal family, and, by Theorem $II.5.3$ of \cite{LehtoVirtanen}, it subconverges in the $C^0$ sense towards an $r$-quasiconformal mapping $\alpha_0:\hat{\Bbb{C}}\rightarrow\hat{\Bbb{C}}$ which preserves $i$, $-i$ and $\Bbb{R}$. This completes the proof.\qed
\proclaim{Proposition \nextprocno}
\noindent $(\alpha_n)^{-1}_\ninn$ subconverges to $\alpha_0$ in the $C^0$ sense.
\endproclaim
\proclabel{PropInversesConverge}
\proof As in the proof of Proposition \procref{PropQuasiConfCompactness}, for all $n\in\Bbb{N}\munion\left\{0\right\}$, we extend $\alpha_n$ to an $r$-quasiconformal homeomorphism of $\hat{\Bbb{C}}$ to itself. Let $K$ and $\Omega$ be subsets of $\hat{\Bbb{C}}$ such that $K$ is compact, $\Omega$ is open, and $\alpha_0^{-1}(K)\subseteq\Omega$. The complements, $\Omega^c$ and $K^c$, are compact and open respectively. Bearing in mind that $\alpha_0$ is bijective:
$$\matrix
&\Omega^c\hfill&\subseteq\alpha_0^{-1}(K)^c\hfill\cr
& &=\alpha_0^{-1}(K^c)\hfill\cr
\Rightarrow\hfill&\alpha_0(\Omega^c)\hfill&\subseteq K^c.\hfill\cr
\endmatrix$$
\noindent Since $(\alpha_n)_\ninn$ converges uniformly to $\alpha_0$, there exists $N\in\Bbb{N}$ such that for $n\geqslant N$:
$$\matrix
&\alpha_n(\Omega^c)\hfill&\subseteq K^c\hfill\cr
\Rightarrow\hfill&\Omega^c\hfill&\subseteq\alpha_n^{-1}(K^c)\hfill\cr
& &=\alpha_n^{-1}(K)^c\hfill\cr
\Rightarrow\hfill&\alpha_n^{-1}(K)\hfill&\subseteq\Omega.\hfill\cr
\endmatrix$$
\noindent The result follows by definition of uniform convergence.\qed
\newhead{Strong Compactness}
\noindent We continue with the notation of the preceeding section. Let $(z_n)_\ninn\in\Bbb{D}$ be a sequence in $\Bbb{D}$ converging to $z_0$ in $\Bbb{D}$. In order to apply Labourie's Theorem, we first show that $(i_n(z_n))_\ninn$ remains in a compact set. We achieve this in two steps:
\proclaim{Proposition \nextprocno}
\noindent There exists a neighbourhood, $U$, of $\varphi_0(z_0)$ in $\Bbb{H}^3\munion\partial_\infty\Bbb{H}^3$ such that, for all $n$:
$$
i_n(z_n)\in U^c.
$$
\endproclaim
\proclabel{PropAvoidsLimitPoint}
\proof For all $n\in\Bbb{N}$, let $\msf{N}_n$ be the normal vector field over $i_n$ compatible with the orientation, and define $I_n:\Bbb{D}\times[0,\infty[\rightarrow\Bbb{H}^3$ by:
$$
I_n(z,t) = \opExp(t\msf{N}_n(z)),
$$
\noindent where $\opExp$ is the exponential map of $\Bbb{H}^3$. Let $B_1$ be the unit ball in $\Bbb{R}^3$. We identify $\Bbb{H}^3$ and $\partial_\infty\Bbb{H}^3$ with $B_1$ and $\partial B_1$ respectively in the canonical manner. For all $n$, $I_n$ extends to a smooth map from $\Bbb{D}\times[0,\infty]$ to $\overline{B}_1$ such that, for all $z\in\Bbb{D}$:
$$
I_n(z,\infty) = \varphi_n(z).
$$
\noindent Let $V$ be a neighbourhood of $z_0$ in $\Bbb{D}$ and suppose that the restriction of $\varphi_0$ to $V$ is a diffeomorphism onto its image. Let $D\subseteq\partial B_1$ be a disc about $\varphi_0(z_0)$ in $\partial B_1$ which is contained in $\varphi_0(V)$. Without loss of generality, we may assume that, for all $n$, $z_n\in V$, the restriction of $\varphi_n$ to $V$ is a diffeomorphism onto its image, and $D\subseteq\varphi_n(V)$.
\medskip
\noindent Let $\Sigma_0$ be the totally geodesic hypersurface in $\Bbb{H}^3$ bounded by $\partial D$. For all $d>0$, let $\Sigma_d$ be the level hypersurface at distance $d$ from $\Sigma$ and let $j_d:\Sigma_d\rightarrow\Bbb{H}^3$ be the canonical embedding. Using elementary hyperbolic geometry, we show that, for all $d$, $\Sigma_d$ has constant Gauss curvature equal to $\opTanh(d)^2$. Denote $d_k=\opTanh^{-1}(\sqrt{k})$. We claim that for all $n$, and for all $d\geqslant d_k$, $\Sigma_d$ lifts to an embedded surface in $\Bbb{D}\times]0,1[$. In other words, there exists $\tilde{\mathj}_{d,n}:\Sigma_d\rightarrow\Bbb{D}\times[0,1]$ such that:
\medskip
\myitem{(i)} $I_n\circ\tilde{\mathj}_{d,n} = j_d$; and
\medskip
\myitem{(ii)} the restriction of $\tilde{\mathj}_{d,n}$ to $\partial_\infty\Sigma_d=\partial D$ coincides with $(\varphi_n|_V)^{-1}$.
\medskip
\noindent Indeed, fix $n$. Let $j_\infty:D\rightarrow\partial_\infty\Bbb{H}^3$ be the canonical embedding. Trivially, $(\Sigma_d,j_d)$ converges to $(D,j_\infty)$ as $d$ tends to $+\infty$. For all $n\in\Bbb{N}$, we define $\tilde{\mathj}_{\infty,n}$ by:
$$
\tilde{\mathj}_{\infty,n} = (\varphi_n|_V)^{-1}\circ j_\infty.
$$
\noindent $\tilde{\mathj}_{\infty,n}$ is trivially a lift of $j_\infty$ satisfying both $(i)$ and $(ii)$. Since $I_n$ is everywhere a local diffeomorphism, for all sufficiently large $d$, $\tilde{\mathj}_{\infty,n}$ perturbs to a lift, $\tilde{\mathj}_{d,n}$ of $j_d$. Bearing in mind that $\hat{\mathi}_n$ is complete, by continuously reducing $d$, $j_d$ can be lifted until we reach some $d=d_0\in[0,\infty]$ where $j_{d_0,n}$ meets $\Bbb{D}\times\left\{0\right\}$ at some point. Suppose that $d_0>d_k$. At the point of contact, $i_n$ is an interior tangent to $j_{d_0}$. However, the Gauss curvature of $\Sigma_d$ is greater than that of $i_n$, which is absurd by the geometric maximum principal. We deduce that $d_0\leqslant d_k$ and the assertion follows.
\medskip
\noindent For sufficiently large $n$, $\varphi_n(z_n)\in D$, and so the geodesic ray in $\Bbb{H}^3$ leaving $i_n(\Bbb{D})$ at $z_n$ in the direction $\msf{N}_n(z_n)$ passes through $\Sigma_{d_k}$. However, this geodesic ray is also the geodesic segment in $\Bbb{H}^3$ joining $i_n(z_n)$ to $\varphi_n(z_n)$. We claim that there exists a neighbourhood, $U$, of $\varphi_0(z_0)$ such that, for sufficiently large $n$, and for any $p\in U$, the geodesic segment in $\Bbb{H}^3$ joining $p$ to $\varphi_n(z_n)$ does not intersect $\Sigma_{d_k}$. Indeed, suppose the contrary. There exists a sequence $(p_n)_\ninn\in\Bbb{H}^3$ converging to $\varphi_0(z_0)$ such that, for all $n$, the geodesic segment in $\Bbb{H}^3$ joining $p_n$ to $\varphi_n(z_n)$ passes through $\Sigma_{d_k}$. However, these geodesic segments converge towards $\varphi_0(z_0)$ in the Haussdorf sense. Thus:
$$
\left\{\varphi_0(z_0)\right\}\minter\Sigma_{d_k}\neq\emptyset.
$$
\noindent This is absurd, and the assertion follows. Thus, for sufficiently large $n$, $i_n(z_n)\notin U$. $U$ is therefore the desired neighbourhood of $\varphi_0(z_0)$, and this completes the proof.\qed
\proclaim{Proposition \nextprocno}
\noindent For any compact subset $K\subseteq\Bbb{D}$, there exists a compact subset $L\subseteq\Bbb{H}^3$ such that, for any $n\in\Bbb{N}$:
$$
i_n(K)\subseteq L.
$$
\endproclaim
\proclabel{PropStaysInCompactSet}
\proof Suppose that contrary. Without loss of generality, there exists $(z_n)_\ninn\in\Bbb{D}$ converging to $z_0\in\Bbb{D}$ such that $(i_n(z_n))_\ninn$ converges to some point $p_0\in\partial_\infty\Bbb{H}^3$. By Proposition \procref{PropAvoidsLimitPoint}, $p_0\neq\varphi_0(z_0)$. Choose $p_0'\in\Bbb{H}^3$, and, for all $n$, let $\Phi_n:\Bbb{H}^3\rightarrow\Bbb{H}^3$ be an isometry such that $\Phi_n(p_n)=p_0'$ and $\Phi_n(\varphi_n(z_n))=\varphi_0(z_0)$.
\medskip
\noindent For all $n$, denote $j_n=\Phi_n\circ\varphi_n$ and let $\hat{\mathj}_{n,r}$ be the Gauss lift of $j_n$ at height $r$. For all $n$, $\hat{\mathj}_{n,r}(z_n)$ lies in the fibre over $p_0'$ which is compact and so, by Theorem \procref{ThmLabCompactness}, there exists a complete, pointed immersed surface $(S_0,\hat{\mathj}_0,q_0)$ towards which $(\Bbb{D},\hat{\mathj}_n,z_n)$ subconverges in the $C^\infty_\oploc$ sense modulo reparametrisation. Let $(\psi_n)_\ninn:S_0\rightarrow\Bbb{D}$ be a sequence of convergence mappings of $(\Bbb{D},\hat{\mathj}_n,z_n)$ with respect to $(S_0,\hat{\mathj}_0,q_0)$. For all $n$, let $J_n$ be the complex structure induced over $\Bbb{D}$ by $\hat{\mathj}_n$. By Proposition \procref{PropConformalRepresentation}, for all $n$, there exists a conformal homeomorphism $\alpha_n:(\Bbb{D},J_n)\rightarrow\Bbb{D}$ sending $z_n$ to $0$. By Lemma \procref{LemmaConvergenceMappingsConverge}, $(\alpha_n\circ\psi_n)_\ninn$ converges in the $C^\infty_\oploc$ sense to a holomorphic mapping $\tilde{\psi}_0:(\Sigma_0,J_0)\rightarrow\Bbb{D}$ sending $q_0$ to $0$. Thus, by Proposition \procref{PropInversesConverge}, $(\psi_n)_\ninn$ converges in the $C^0_\oploc$ sense to $\psi_0:=\alpha_0^{-1}\circ\tilde{\psi}_0$. Trivially $\varphi_0(q_0)=z_0$.
\medskip
\noindent Since $p_0\neq\varphi_0(z_0)$, $(\Phi_n)_\ninn$ converges locally uniformly over $(\Bbb{H}^3\munion\partial_\infty\Bbb{H}^3)\setminus\left\{p_0\right\}$ to the constant mapping sending every point to $\varphi_0(z_0)$. There thus exists a neighbourhood, $U$, of $z_0$ over which $(\Phi_n\circ\varphi_n)_\ninn$ converges to a constant map, and so there exists a neighbourhood, $V$, of $q_0$ over which $(\Phi_n\circ\varphi_n\circ\psi_n)_\ninn$ also converges to a constant map. However, for all $n$, bearing in mind that $\overrightarrow{n}$ commutes with $\Phi_n$:
$$\matrix
\Phi_n\circ\varphi_n\circ\psi_n\hfill&=\Phi_n\circ\overrightarrow{n}\circ\hat{\mathi}_n\circ\psi_n\hfill\cr
&=\overrightarrow{n}\circ\Phi_n\circ\hat{\mathi}_n\circ\psi_n\hfill\cr
&=\overrightarrow{n}\circ\hat{\mathj}_n\circ\psi_n.\hfill\cr
\endmatrix$$
\noindent This converges in the $C^\infty_\oploc$ sense to $\overrightarrow{n}\circ\hat{\mathj}_0$ which is therefore constant over $V$. This is absurd, since, in both cases given in Theorem \procref{ThmLabCompactness}, $\overrightarrow{n}\circ\hat{\mathj}_0$ is a local homeomorphism. The result follows.\qed
\medskip
\noindent We are now in a position to prove:
\proclaim{Proposition \nextprocno}
\noindent There exists a complete immersion $\hat{\mathi}_0:\Bbb{D}\rightarrow S_r\Bbb{H}^3$ towards which $(\hat{\mathi}_n)_\ninn$ subconverges in the $C^\infty_\oploc$ sense. Moreover, $i_0=\pi\circ\hat{\mathi}_0$ is an immersion.
\endproclaim
\proclabel{PropStrongCompactness}
\proof For all $n$, let $A_n$ be the shape operator of $i_n$. We claim that for every compact subset $K\subseteq\Bbb{D}$, there exists $B>0$ such that for all $z\in K$, and for all $n$, $\|A_n(z)\|\leqslant B$. Indeed, suppose the contrary. Then, without loss of generality, there exists a sequence $(z_n)_\ninn\in\Bbb{D}$ converging to a limit $z_0\in\Bbb{D}$ such that $(\|A_n(z)\|)_\ninn\rightarrow+\infty$. By Proposition \procref{PropStaysInCompactSet}, there exists a compact subset, $L\subseteq\Bbb{H}^3$ such that $i_n(z_n)\in L$ for all $n$. Thus, by Theorem \procref{ThmLabCompactness}, there exists a complete, pointed immersion $(S_0,\hat{\mathi}_0,p_0)$ towards which $(\Bbb{D},\hat{\mathi}_n,z_n)$ subconverges in the $C^\infty_\oploc$ sense modulo reparametrisation. Since $(\|A_n\|)_\ninn\rightarrow+\infty$, $\pi\circ\hat{\mathi}_0$ cannot be an immersion at $p_0$, and the immersed surface $(S_0,\hat{\mathi}_0)$ is therefore tubular.
\medskip
\noindent Let $(\psi_n)_\ninn:S_0\rightarrow\Bbb{D}$ be a sequence of convergence mappings of $(\Bbb{D},\hat{\mathi}_n,z_n)_\ninn$ with respect to $(S_0,\hat{\mathi}_0,p_0)$. For all $n$, let $J_n$ be the complex structure induced over $\Bbb{D}$ by $\hat{\mathi}_n$. By Proposition \procref{PropConformalRepresentation}, for all $n$, there exists a conformal homoemorphism $\alpha_n:(\Bbb{D},J_n)\rightarrow\Bbb{D}$ sending $z_n$ to $0$. By Lemma \procref{LemmaConvergenceMappingsConverge}, $(\alpha_n\circ\psi_n)_\ninn$ subconverges in the $C^\infty_\oploc$ sense to a holomorphic mapping $\tilde{\psi}_0:(S_0,J_0)\rightarrow\Bbb{D}$ sending $q_0$ to $0$. Since $(S_0,\hat{\mathi}_0)$ is tubular, in particular, it is of parabolic type, and so $\tilde{\psi}_0$ is constant. However, by Proposition \procref{PropInversesConverge}, $(\psi_n)_\ninn$ converges in the $C^0_\oploc$ sense to $\psi_0:=\alpha_0^{-1}\circ\tilde{\psi}_0$. Since $\tilde{\psi}_0$ is constant, so is $\psi_0$. Thus, in particular, $(\varphi_n\circ\psi_n)_\ninn$ converges in the $C^0_\oploc$ sense to a constant mapping. However, for all $n$:
$$
\varphi_n\circ\psi_n = \overrightarrow{n}\circ\hat{\mathi}_n\circ\psi_n.
$$
\noindent Since $(\hat{\mathi}_n\circ\psi_n)_\ninn$ converges in the $C^\infty_\oploc$ sense to $\hat{\mathi}_0$, we conclude that $\overrightarrow{n}\circ\hat{\mathi}_0$ is a constant mapping. This is absurd and the assertion follows.
\medskip
\noindent By Theorem \procref{ThmLabCompactness}, there exists a complete, pointed immersion $(S_0,\hat{\mathi}_0,p_0)$ towards which $(\Bbb{D},\hat{\mathi}_n,0)$ subconverges in the $C^\infty$ sense modulo reparametrisation. Let $(\psi_n)_\ninn:S_0\rightarrow\Bbb{D}$ be a sequence of convergence mappings of $(\Bbb{D},\hat{\mathi}_n,0)$ with respect to $(S_0,\hat{\mathi}_0,p_0)$. As before $(\psi_n)_\ninn:S_0\rightarrow\Bbb{D}$ converges in the $C^0_\oploc$ sense to some $\psi_0:S_0\rightarrow\Bbb{D}$. Taking limits, we obtain:
$$
\varphi_0\circ\psi_0 = \overrightarrow{n}\circ\hat{\mathi}_0,
$$
\noindent and so $\varphi_0\circ\psi_0$ is everywhere a local diffeomorphism. It follows that $\psi_0$ is a local diffeomorphism and $(\psi_n)_\ninn$ subconverges to $\psi_0$ in the $C^\infty_\oploc$ sense.
\medskip
\noindent We claim that $\psi_0$ is a covering map. Let $\gamma:[0,1]\rightarrow\Bbb{D}$ be a curve such that $\gamma(0)=0$. Let $\hat{g}$ be the Sasaki metric over $T\Bbb{H}^3$. For all $n$, let $l_n(\gamma)$ be the length of $\gamma$ with respect to $\hat{\mathi}_n^*\hat{g}$. By Proposition \procref{PropStaysInCompactSet}, there exists a compact subset, $L\subseteq\Bbb{H}^3$ such that $i_n(\gamma([0,1]))\subseteq L$ for all $n$. Since, in addition, $(\varphi_n)_\ninn$ converges and $(A_n)_\ninn$ is uniformly bounded, we show that $(\hat{\mathi}_n\circ\gamma)_\ninn$ is uniformly bilipschitz. There thus exists $R>0$ such that, for all $n$:
$$
l_n(\gamma)\leqslant R.
$$
\noindent Thus, for sufficiently large $n$, $\gamma$ lifts through $\psi_n$ to a smooth curve $\gamma_n\subseteq S_0$ of length at most $2R$. Taking limits, we see that $\gamma$ lifts through $\psi_0$ to a smooth curve in $S_0$, and the assertion follows.
\medskip
\noindent Since $\psi_0$ is a covering map and $\Bbb{D}$ is simply connected, $\psi_0$ is a diffeomorphism. Let $K\subseteq\Bbb{D}$ be compact. Let $\Omega\subseteq\Bbb{D}$ be a relatively compact neighbourhood of $\Bbb{D}$. There exists $N\in\Bbb{N}$ such that, for all $n\geqslant N$, $\psi_n$ has a smooth inverse, $\beta_n:\Omega\rightarrow S_0$. Moreover, $(\beta_n)_\ninn$ converges in the $C^\infty_\oploc$ sense to $(\psi_0|_\Omega)^{-1}$. For all $n$:
$$
\hat{\mathi}_n = \hat{\mathi}_n\circ\psi_n\circ\beta_n.
$$
\noindent Thus $(\hat{\mathi}_n)_\ninn$ converges in the $C^\infty_\oploc$ sense to $\hat{\mathi}_0\circ(\psi_0|_\Omega)^{-1}$. This completes the proof.\qed
\medskip
\noindent Theorems \procref{ThmExistence} and \procref{ThmContinuousDependance} follow readily:
\medskip
{\bf\noindent Proof of Theorem \procref{ThmExistence}:~}For all $t<1$, define $\alpha_t:\Bbb{D}\rightarrow\Bbb{D}$ by:
$$
\alpha_t(z) = tz.
$$
\noindent By Theorem $E$ of \cite{Labourie}, for all $t$, there exists $i_t:\Bbb{D}\rightarrow\Bbb{H}^3$ such that:
\medskip
\myitem{(i)} $i_t$ has constant Gaussian curvature equal to $k$;
\medskip
\myitem{(ii)} $\hat{\mathi}_{t,r}$ is complete; and
\medskip
\myitem{(iii)} $\overrightarrow{n}\circ\hat{\mathi}_{t,r}=\varphi\circ\alpha_t$.
\medskip
\noindent Letting $t$ converge to $1$, existence follows by Theorem \procref{ThmCompactness}. Uniqueness follows by Theorem $A$ of \cite{Labourie}. This completes the proof.\qed
\medskip
{\bf\noindent Proof of Theorem \procref{ThmContinuousDependance}:~}Since $\varphi_n$ is conformal for all $n$, $(\varphi_n)_\ninn$ converges to $\varphi_0$ in the $C^\infty_\oploc$ sense. Let $(i_{k_n})_\ninn$ be a subsequence of $(i_n)_\ninn$. By Theorem \procref{ThmCompactness}, there exists $i'_0$ towards which $(i_{k_n})_\ninn$ subconverges in the $C^\infty_\oploc$ sense such that:
\medskip
\myitem{(i)} $i_0'$ has constant Gaussian curvature equal to $k$;
\medskip
\myitem{(ii)} $(\hat{\mathi}'_{0,r})$ is complete; and
\medskip
\myitem{(iii)} $\overrightarrow{n}\circ\hat{\mathi}_{0,r}'=\varphi_0$.
\medskip
\noindent By uniqueness (c.f. Theorem $A$ of \cite{Labourie}), $i_0'=i$. Thus, every subsequence of $(i_n)_\ninn$ has a subsubsequence converging to $i_0$, and so $(i_n)_\ninn$ itself converges to $i_0$. This completes the proof.\qed
\goodbreak
\newhead{Bibliography}
{\leftskip = 5ex \parindent = -5ex
\leavevmode\hbox to 4ex{\hfil \cite{BallSchGrom}}\hskip 1ex{Ballmann W., Gromov M., Schroeder V., {\sl Manifolds of nonpositive curvature}, Progress in Mathematics, 61, Birkh\"auser Boston Inc., Boston, MA, (1985)}
\medskip
\leavevmode\hbox to 4ex{\hfil \cite{Brezis}}\hskip 1ex{Brezis H., {\sl Functional analysis, Sobolev spaces and partial differential equations}, Universitext, Springer, New York, (2011)}
\medskip
\leavevmode\hbox to 4ex{\hfil \cite{CaffNirSpruckV}}\hskip 1ex{Caffarelli L., Nirenberg L., Spruck J., Nonlinear second-order elliptic equations. V. The Dirichlet problem for Weingarten hypersurfaces, {\sl Comm. Pure Appl. Math.} {\bf 41} (1988), no. 1, 47--70}
\medskip
\leavevmode\hbox to 4ex{\hfil \cite{Gauss}}\hskip 1ex{Gauss C. F., Disquisitiones generales circa superficies curvas, Oct. 8 1827}
\medskip
\leavevmode\hbox to 4ex{\hfil \cite{Gromov}}\hskip 1ex{Gromov M., Pseudoholomorphic curves in symplectic manifolds, {\sl Invent. Math.} {\bf 82} (1985), no. 2, 307--347}
\medskip
\leavevmode\hbox to 4ex{\hfil \cite{GuanSpruck}}\hskip 1ex{Guan B., Spruck J., The existence of hypersurfaces of constant Gauss curvature with prescribed boundary, {\sl J. Differential Geom.} {\bf 62} (2002), no. 2, 259--287}
\medskip
\leavevmode\hbox to 4ex{\hfil \cite{ImaTani}}\hskip 1ex{Imayoshi Y., Taniguchi M., {\sl An introduction to Teichm\"uller space}, Springer-Verlag, Tokyo, (1992)}
\medskip
\leavevmode\hbox to 4ex{\hfil \cite{LabourieB}}\hskip 1ex{Labourie F., Probl\`eme de Minkowski et surfaces \`a courbure constante dans les vari\'et\'es hyperboliques, {\sl Bull. Soc. Math. France} {\bf 119} (1991), no. 3, 307--325}
\medskip
\leavevmode\hbox to 4ex{\hfil \cite{Labourie}}\hskip 1ex{Labourie F., Un lemma de Morse pour les surfaces convexes, {\sl Invent. Math.} {\bf 141} (2000), no. 2, 239--297}
\medskip
\leavevmode\hbox to 4ex{\hfil \cite{LehtoVirtanen}}\hskip 1ex{Lehto O., Virtanen K. I., {\sl Quasiconformal mappings in the plane}, die Grundlehren der mathematischen Wissenschaften, {\bf 126}, Springer-Verlag, New-York, Heidelberg, (1973)}
\medskip
\leavevmode\hbox to 4ex{\hfil \cite{Martinez}}\hskip 1ex{Corro A. V., Mart\'inez A., Mil\'an F., Complete flat surfaces with two isolated singularities in hyperbolic 3-space, {\sl J. Math. Anal. Appl.} {\bf 366} (2010), no. 2, 582--592}
\medskip
\leavevmode\hbox to 4ex{\hfil \cite{Pogorelov}}\hskip 1ex{Pogorelov A. V., {\sl Extrinsic geometry of convex surfaces}, Translations of Mathematical Monographs, Vol. 35. American Mathematical Society, Providence, R.I., (1973)}
\medskip
\leavevmode\hbox to 4ex{\hfil \cite{RosSpruck}}\hskip 1ex{Rosenberg H., Spruck J., On the existence of convex hypersurfaces of constant Gauss curvature in hyperbolic space, {\sl J. Differential Geom.} {\bf 40} (1994), no. 2, 379--409}
\medskip
\leavevmode\hbox to 4ex{\hfil \cite{Sasaki}}\hskip 1ex{Sasaki S., On the differential geometry of tangent bundles of Riemannian manifolds, {\sl T\^ohoku Math. J. (2)} {\bf 10} (1958), 338--354}
\medskip
\leavevmode\hbox to 4ex{\hfil \cite{SmiPKS}}\hskip 1ex{Smith G., Pointed k-surfaces, {\sl Bull. Soc. Math. France} {\bf 134}, no. 4, (2006), 509--557}
\medskip
\leavevmode\hbox to 4ex{\hfil \cite{SmiFCS}}\hskip 1ex{Smith G., Moduli of Flat Conformal Structures of Hyperbolic Type, to appear in {\sl Geom. Dedicata}}
\medskip
\leavevmode\hbox to 4ex{\hfil \cite{SmiPPH}}\hskip 1ex{Smith G., Compactness for immersions of prescribed Gaussian curvature II - geometric aspects, arXiv:1002.2982}
\medskip
\leavevmode\hbox to 4ex{\hfil \cite{TrudWang}}\hskip 1ex{Trudinger N. S., Wang  X. J., On locally convex hypersurfaces with boundary, {\sl J. Reine Angew. Math.} {\bf 551} (2002), 11--32}
\par}%
\enddocument